# A Multi-Period Market Design for Markets with Intertemporal Constraints

Jinye Zhao, Member, *IEEE*, Tongxin Zheng, Senior Member, *IEEE*, and Eugene Litvinov, *Fellow, IEEE*

*Abstract*—The participation of renewable, energy storage, and resources with limited fuel inventory in electricity markets has created the need for optimal scheduling and pricing across multiple market intervals for resources with intertemporal constraints. In this paper, a new multi-period market model is proposed to enhance the efficiency of markets with such type of resources. It is also the first market design that links a forward market and a spot market through the coordination of schedule and price under the multi-period paradigm, achieving reliability, economic efficiency and dispatch-following incentives simultaneously. The forward market solves a multi-period model with a long look-ahead time horizon whereas the spot market solves a series of multi-period dispatch and pricing problems with a shorter look-ahead time horizon on a rolling basis. By using the forward schedules and opportunity costs of intertemporal constraints as a guideline, the spot market model is able to produce economically efficient dispatch solutions as well as prices that incentivize dispatch following under the perfect forecast condition. The proposed scheme is applied to the dispatch and pricing of energy storage resources. Numerical experiments show that the proposed scheme outperforms the traditional myopic method in terms of economic efficiency, dispatch following and reliability.

## I. INTRODUCTION

INTERTEMPROAL constraints are inherent to almost all the resources participating in electricity markets today. For instance, every generating unit has a ramp rate, which limits the output change within a certain period. Energy storage resources have limited storage capability, limiting the amount of charge and discharge in each operating cycle. Units such as oil, pipeline constrained gas-fired units and pondage hydro may have limited fuel inventories and/or emission limits, which prevent them from operating at their full capacity for the entire operating day.

Intertemporal constraints couple the markets in different time intervals. The market outcome in one market interval may affect another. Under several recent industry trends in the U.S., such temporal market coupling has become tighter:

- The "duck curve" load shape resulting from a large amount of PV installation [1] requires more ramping capacity from conventional units to cope with a deep midday drop in the net load. Consequently, market clearing are more frequently constrained by the intertemporal linkages among different market intervals, making it difficult for market participants to predict these coupling intervals and internalize the corresponding opportunity cost into their bids.
- The states' effort to pursue their renewable portfolio standards helps spur a large amount of energy storage resources (ESR) [2]–[4]. The emerging storage technologies, such as batteries and flywheels, are capable of running many charge/discharge cycles per day. However, the outcomes of self-managed ESRs may not be efficient and reliable due to the dynamics in the real-time market. An alternative would be ISO-managed ESRs, whose intertemporal constraints are explicitly modeled in the ISO scheduling and pricing processes.
- As the nation increasingly relies on the natural-gas-fired units that are largely connected to the gas pipelines rather than having local gas storage, disruptions on the gas pipeline systems can have a significant impact on the reliable operation of the power system [5]. The value of the fuel storage becomes more prominent, and optimal scheduling of these resources across different market intervals is critical.

The above industry trends highlight the implication of intertemporal linkages for optimal operation of power systems and call for careful studies of scheduling and pricing methods or market design for markets with intertemporal constraints.

## II. PRIOR WORK AND CONTRIBUTIONS

All the ISOs in the U.S. run a day-ahead (DAM) and real-time market (RTM). The DAM, a forward market, typically contains 24 market periods. The schedules and prices for all the DAM periods are simultaneously produced by a 24-hour market clearing model, respecting intertemporal constraints. As a result, the DAM clearing prices reflect intertemporal opportunity cost, which is the opportunity cost a resource incurred due to the fact that the production in one period affects another via intertemporal linkages. On the other hand, the RTM, the spot market, is typically a single-period market, whose clearing prices may not reflect intertemporal opportunity costs. RT schedules are determined by optimizing either a single- or multi-period problem depending on the ISOs' tariffs. The existing RTM scheduling and pricing models are reviewed in Subsection II.A.

Under the existing market designs in the U.S., the coordination between the DAM and RTM is weak in terms of dispatch and pricing. The DAM schedules are not considered in the RTM scheduling process. By relying on the information only in a short RT time frame, the RTM schedules may impair reliability and efficiency, especially for inflexible generators, like oil or coal units, which require several hours to reach their full capacity as well as for storage hydro and limited energy resources, whose schedules need to be managed over the day. Furthermore, the existing RTM prices can be inconsistent with the DAM prices. For example, unlike the DAM prices derived from a 24-hour clearing model, the RTM price in a peak hour may not account for the intertemporal opportunity costs of ramping up a unit out of merit in the previous hours to meet the peak demand. Such RTM price can discourage dispatch-following incentives[1] in the short term and undervalue flexible resources used to meet the peak load in the long term as pointed out in [6].[2]

Although there exists a large body of literature on the multi-period scheduling models, there are few publications that examine the scheduling coordination between a forward market and the RTM, the linkage between forward market prices and RTM prices, as well as the impact of RTM pricing/settlement schemes on dispatch following incentives.

*A. Prior work*

   *1) The myopic approach*

The myopic approach in the RTM solves a single-period economic dispatch problem that minimizes the cost of meeting load over a 5- or 10-minute period, without considering the cost or the system's ability of meeting demand in future dispatch intervals. As a result, system operators may have to take manual actions to adjust dispatch solutions to accommodate expected future conditions. Such manual actions are subjective and can be suboptimal or even infeasible. In addition, market clearing prices generated based on a single-period model are often inconsistent with manual actions, creating disincentives for units to follow their dispatch instructions, and subsequently increasing the need for out-of-market payments.

   *2) The multi-period approach*

Multi-period market clearing models have been drawing much attention because of their capability of optimizing the system with a set of intertemporal constraints across several time intervals. A few ISOs in the U.S. have implemented multi-period dispatch in the RTM [7]–[9].

   *i) Challenges of existing multi-period scheduling models*

Various multi-period scheduling models have been extensively studied with the objectives of improving dispatch efficiency and reliability, e.g. [10]–[13], etc. However, when the market coupling involves many market intervals, the sizes of these models become too large to be solved quickly to provide solutions every five minutes in the RTM. In practice, in order to reduce the model size, the look-ahead horizon of these models is typically set at one hour, as in New York and Californian ISOs, and the time resolution is increased to reduce the number of dispatch periods in the multi-period scheduling model at the expense of operational efficiency and reliability.

   *ii) Challenges of existing multi-period pricing models*

Current market designs do not have efficient pricing or settlement mechanisms for a resource to recover its intertemporal opportunity costs via market clearing prices. Consequently, resources are compelled to internalize their opportunity costs into bids and offers, but estimating opportunity costs is a challenging task as pointed out in [14]. Therefore, it is important to develop an efficient pricing scheme that properly reflects the impact of intertemporal constraints in the market clearing prices. The existing multi-period pricing work can be categorized into three types, namely the first-period only settlement, the multi-settlement, and the prior-cost based pricing. These pricing approaches all focus on the RTM, but do not exploit the connection between the DAM and RTM.

- Under the first-period only settlement, a multi-period pricing problem is run at each market clearing. Only the first period is settled, and the prices for the later periods are advisory and not financially binding [15]. This approach has been used in New York ISO [7] and California ISO [8]. Because opportunity costs embedded in the advisory intervals are never settled upon, as noted in [16], the market participants may have incentives to deviate from dispatch instructions, decreasing operational efficiency.

- Unlike the first-period only settlement approach, all the periods are financially binding under the multi-settlement approach [17]. However, because the multi-settlement approach considers the past dispatch decisions irrelevant to price determination, the resulting prices do not necessarily accurately reflect the total cost of serving the load. Consequently, this approach can distort dispatch-following incentives in the short term as shown in [16], as well as fail to send right investment signals in the long run.

- The prior-cost based pricing method in [18] and [19] incorporates the costs associated with the past dispatch decisions, and then determines prices by re-optimizing the social welfare for the entire look-ahead time horizon [18] or the remainder of the time horizon [19]. Only the first interval is settled. This approach is shown to alleviate the dispatch-following incentive issue in RT. However, if the RTM look-ahead time horizon is long, e.g. several hours or a day, the prior-cost based approach would be impractical because it is very computationally challenging to solve such large-scale problems in the RT operation.

*B. Contributions*

To address the issue with the weak coordination between the forward market and RTM under the existing approaches, we propose a coupled multi-period market design that strengthens the coordination between these two markets. Under the proposed design, the forward market solves a multi-period problem with a long look-ahead horizon as the existing forward markets in the U.S. In the RTM, the RT prices are determined by solving linear programming problems with a shorter look-ahead horizon. The RT pricing problems coordinate the RT prices with forward prices to better reflect intertemporal opportunity costs, and the resulting RT prices provide better dispatch-following incentives. In addition, a new RT scheduling problem is employed to take into account the schedules obtained in the forward market. It ensures efficient and reliable system operation under the perfect forecast condition as well as is computationally efficient for real-world implementations. A multi-settlement system is adopted for settling each market clearing, reducing market participants' risk exposure. The mechanism of this multi-settlement system is similar to the multi-settlement pricing approach proposed in [17]. However, unlike the multi-settlement pricing approach, the proposed clearing prices more accurately reflect the load serving costs, which in turn induce dispatch-following incentives.

---

[1] A resource has dispatch-following incentives under a given clearing price if the resource's profit would not be better off by not following the ISO's dispatch instructions under the clearing price.

[2] Although unit commitment processes are run periodically in the RTM to commit additional units or pre-ramp resources for the reliability purpose, there are no market prices associated with these commitment and dispatch instructions. Therefore, the existing real-time unit commitment processes do not properly provide dispatch-following incentives.

Although the proposed approach is seemingly comparable to bi-level or multi-stage stochastic models, it is fundamentally different. In bi-level approach [20]–[22], no pricing or scheduling mechanisms are offered while a three-stage stochastic model in [23] is used only for DAM clearing, but does not consider intertemporal coupling in RTM.

The proposed pricing design can be distinguished from [19] from two perspectives. First, [19] investigates the RTM pricing model in isolation whereas the proposed design considers both forward market and RTM, leading to improvement in efficiency and reliability. Second, the proposed pricing design can be considered as an evolution of [19] for its flexibility to cover a spectrum of deployment time frames without sacrificing computational efficiency, economic efficiency and dispatch-following incentives, making it more practical and general than [19].

The specific contributions of this paper are the following:
a. We develop a practical multi-period market clearing model for scheduling and pricing of electricity markets with intertemporal constraints. Under the proposed approach, the ISO is able to maintain operational reliability and market efficiency without sacrificing computational efficiency.
b. This is the first paper that considers the coordination between a forward market and RTM in terms of scheduling and pricing under the multi-period paradigm. The incorporation of intertemporal opportunity costs in the RTM pricing scheme makes the RT prices more closely reflect the load serving cost, maintaining the pricing consistency between the forward prices and RT prices. It also reduces the need of participants incorporating the opportunity cost in their bids and offers.
c. The proposed pricing scheme provides proper dispatch-following incentives, and the proposed scheduling scheme produces economically efficient dispatch results under the perfect forecast conditions. Other properties of the proposed scheme are also studied under the imperfect forecast conditions.
d. We apply the proposed scheme to a stylized small system with ESRs as well as real-world ISO New England (ISO NE) system. The advantages of the proposed scheme in terms of economic efficiency, reliability and dispatch-following incentives are demonstrated in the numerical studies.

The rest of the paper is organized as follows. The proposed coordinated market clearing model as well as its properties are presented in Section III. The numerical studies are described in Section IV. The paper is concluded in Section V.

## III. THE PROPOSED COORDINATED MULTI-PERIOD MARKET CLEARING FRAMEWORK

Under the proposed scheme, the forward market takes place at $t_0$ as illustrated in Fig. 1. A multi-period model is solved to meet the forecasted system condition at the entire planning time horizon, and simultaneously generates the forward clearing quantities and prices for each time period. It also obtains intertemporal opportunity costs for each resource under the forecasted system condition.

The RTM has a shorter look-ahead time horizon, which is called *sub time horizon*, and clears on a rolling basis. Each market clearing results in a set of cleared quantities together with a set of corresponding prices. The cleared quantities and prices are settled after each market clearing under a multi-settlement system. The shorter time horizon in the RTM is able to reduce computational burden as well as simplify the settlement process. In order to compensate the narrower vision of the RT market clearing, the forward schedules are utilized to guide the RT scheduling model. The intertemporal opportunity costs are incorporated into the RT pricing model to maintain the consistency of prices between the forward market and RTM, providing much needed dispatch-following incentives.

Even though the forecast in the forward market may not be perfect, it can still be argued that the forward solution is the best possible solution given system operators' beliefs. The RT scheduling and pricing problems are designed to follow the forward scheduling and price trajectories, and at the same time it is left with the freedom to adjust schedules and prices the best way it can in the sub time horizon given the newly revealed information.

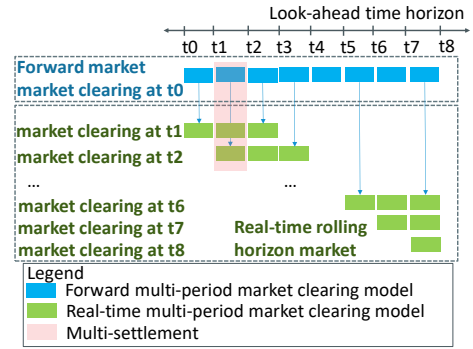

**Figure 1: The proposed multi-period market clearing framework**

### A. The forward market

#### 1) The forward market clearing model

The forward multi-period market clearing model can be cast as the following linear programming problem[3]:

---

[3] We focus on linear programming problems in this paper. The integer decisions and other non-convex constraints are not considered because these features are closely related to the non-convexity pricing issue, which is not the focus of this paper. In the presence of the non-convexities, there may be no set of prices that supports the dispatch to achieve a competitive partial equilibrium. In the current practice, one solution to make the dispatch consistent with the incentives in the market is to make up any losses for dispatched resources via the uplift payments. In addition, some ISOs have modified the traditional marginal pricing scheme to incorporate non-convex costs such as start-up and no-load costs into prices and redefined uplift payments. Although the proposed pricing scheme addresses the issues associated with intertemporal constraints, it can be incorporated into the existing practice of dealing with non-convexities.

$$\min_{x_1,\ldots,x_T} \sum_{t\in\mathcal{T}} c_t^\mathsf{T} x_t \quad (1)$$

$$\text{s.t.} \quad B_t x_{t-1} + A_t x_t \le b_t \quad t\in\mathcal{T},\ (\pi_{t-1:t}) \quad (2)$$

$$G_t x_t \le g_t \quad t\in\mathcal{T},\ (\lambda_t) \quad (3)$$

$$F_t x_t \le f_t \quad t\in\mathcal{T},\ (\mu_t) \quad (4)$$

The objective of the multi-period problem (1) is to minimize the overall production cost across the entire time horizon $t\in\mathcal{T}=\{1,\ldots,T\}$. Constraint (2) represents the intertemporal constraints. Constraint (3) is the system-level time-independent constraint, representing the energy balance, reserve requirement and transmission constraints. Constraint (4) represents time-independent constraints at the individual resource level, e.g. resources' operational capacity.

We use the superscript $u$ to denote the optimal primal and dual solutions to the forward problem (1)-(4): $x^u := \{x_1^u, \ldots, x_T^u\}$ and $\Lambda^u := (\pi^u, \lambda^u, \mu^u)$. The forward clearing quantity is $x^u$, and the market clearing price is $LMP_t^u := G_t^\mathsf{T} \lambda_t^u$ for $t\in\mathcal{T}$.

*2) Properties of the forward market clearing model*

In this subsection, we discuss the properties of the forward market clearing results. Consider the profit maximization problem of a price-taking resource $i$ under the clearing price $LMP^u$:

$$\max_{x_i} \sum_{t\in\mathcal{T}} (LMP_t^u - c_{i,t}) x_{i,t} \quad (5)$$

$$\text{s.t.} \quad B_{i,t} x_{i,t-1} + A_{i,t} x_{i,t} \le b_{i,t} \quad t\in\mathcal{T},\ (\pi_{i,t-1:t}) \quad (6)$$

$$F_{i,t} x_{i,t} \le f_{i,t} \quad t\in\mathcal{T},\ (\mu_{i,t}) \quad (7)$$

Based on the definition in [24], the clearing quantity $x^u$ and price $LMP^u$ constitute a *competitive equilibrium* if a) for each resource $i$, $x_i^u$ solves the profit maximization problem (5)-(7), and b) supply and demand are balanced. Property 1 indicates that the forward clearing results achieve a competitive equilibrium.

**Property 1:** Assuming that all the resources bid truthfully, the forward clearing quantity $x^u$ and price $LMP^u$ constitute a competitive equilibrium.

Property 1 can be easily proved by matching the Karush-Kuhn-Tucker (KKT) conditions of the upper-level overall problem (1)-(4) with the KKT conditions of each individual profit maximization problem (5)-(7) as well as the market clearing condition (3). We omit the details here.

Next, consider the Lagrangian function of the profit maximization problem (5)-(7):

$$L_i(x_i,\pi_i,\mu_i) = \sum_{t\in\mathcal{T}} \begin{bmatrix} (LMP_t^u - c_{i,t})\, x_{i,t} + \\ \pi_{i,t-1:t}^\mathsf{T}(B_{i,t} x_{i,t-1} + A_{i,t} x_{i,t} - b_{i,t}) \\ + \mu_{i,t}^\mathsf{T}(F_{i,t} x_{i,t} - f_{i,t}) \end{bmatrix} \quad (8)$$

It is easy to show that the forward dual solution ($\pi_i^u, \mu_i^u$) is the Lagrangian multipliers for (8). Therefore, by the first-order condition, we have

$$LMP_t^u - c_{i,t} + A_{i,t}^\mathsf{T} \pi_{i,t-1:t}^u + B_{i,t+1}^\mathsf{T} \pi_{i,t:t+1}^u + F_{i,t}^\mathsf{T} \mu_{i,t}^u = 0 \quad (9)$$

If $\mu_{i,t}^u$ is zero, it means that resource $i$'s output is between its minimum and maximum capacities. Thus, the resource is a marginal resource. As a result, equation (9) can be written as follows:

$$LMP_t^u = c_{i,t} - A_{i,t}^\mathsf{T} \pi_{i,t-1:t}^u - B_{i,t+1}^\mathsf{T} \pi_{i,t:t+1}^u. \quad (10)$$

Notice that $LMP_t^u$ and $c_{i,t}$ are no longer equal to each other in the time period $t$ as is expected for a marginal resource in a single-period problem. This is because the cost $c_{i,t}$ is not the true marginal cost – it is only the marginal production cost, e.g. the fuel cost. The true marginal cost is the right hand side of (10), including the marginal production cost $c_{i,t}$ as well as the intertemporal opportunity cost $-B_{i,t+1}^\mathsf{T} \pi_{i,t:t+1}^u - A_{i,t}^\mathsf{T} \pi_{i,t-1:t}^u$.

In particular, the *intertemporal opportunity cost* in the time period $t$ is a resource's overall foregone profit if it chooses one more unit production in $t$, which in turn affects the production in other time periods due to intertemporal constraints. The *marginal profit* in the time period $t$ is the additional profit in the time period $t$ of choosing one more unit of production in that time period. The marginal profit can be either positive (marginal gain) or negative (marginal loss). Equation (10) implies the following property.

**Property 2:** Assuming a resource is marginal in the time period $t$, the marginal profit $LMP_t^u - c_{i,t}$ compensates the intertemporal opportunity cost $-B_{i,t+1}^\mathsf{T} \pi_{i,t:t+1}^u - A_{i,t}^\mathsf{T} \pi_{i,t-1:t}^u$ for this resource.

Proof: The marginal profit in the time period $t$ is $LMP_t^u - c_{i,t}$ for resource $i$ by the definition. In addition, because $\pi_i^u$ is the Lagrangian multiplier for the profit maximization problem (5)-(7) based on Property 1, it follows by sensitivity analysis in [25] that the overall profit change, i.e. the objective value change in (5), would be $-B_{i,t+1}^\mathsf{T} \pi_{i,t:t+1}^u - A_{i,t}^\mathsf{T} \pi_{i,t-1:t}^u$ if resource $i$ were allowed to provide one additional unit from $x_{i,t}^u$ at time period $t$. Hence, $-B_{i,t+1}^\mathsf{T} \pi_{i,t:t+1}^u - A_{i,t}^\mathsf{T} \pi_{i,t-1:t}^u$ represents the intertemporal opportunity cost. Together with equation (10), it implies that the marginal profit of the time period $t$ equals the intertemporal opportunity cost. Q.E.D.

*B. The RTM*

*1) The RTM clearing model*

The RTM is a rolling-horizon market. Each RT market clearing looks ahead a shorter time horizon rather than the full horizon of the operating day. We use the subscript $s$ to denote a RTM sub time horizon $\mathcal{S} = \{t_s^0, \ldots, t_s^*\}$ starting at $t_s^0$ and ending at $t_s^*$. The subscript $s'$ denotes the past time horizon $\mathcal{S}' = \{1,...,t_s^0 - 1\}$, and the subscript $s''$ denotes the future time horizon $\mathcal{S}'' = \{t_s^* + 1,...,T\}$. Without loss of generality, it is assumed that the parameter $g$ is uncertain. The realized uncertainty in the past is denoted by $g_{s'}^*$; new information is revealed for the current sub time horizon, resulting in a new forecast $\hat{g}_s$; and the future forecast is assumed to remain unchanged as $g_{s''}$.

Different from the forward market clearing, the RTM clearing involves two processes, namely, a scheduling process, where the market clearing quantity is determined to maximize the social welfare in the sub time horizon, and a pricing process, where the market clearing price is determined for each time period in the sub time horizon.

The proposed RT scheduling problem at the sub time horizon $\mathcal{S}$ is defined by (11), using the realized dispatch $x_{t_s^0-1}^*$ and the forward schedule $x_{t_s^*+1}^u$ of the subsequent time period as the boundary conditions:

$$SP_s(x_{t_s^0-1}^*, x_{t_s^*+1}^u) := \min_{x_s} c_s^\mathsf{T} x_s$$
$$\text{s.t. } B_{t_s^0} x_{t_s^0-1}^* + A_{t_s^0} x_{t_s^0} \le b_{t_s^0}$$
$$B_t x_{t-1} + A_t x_t \le b_t, \quad t = t_s^0 + 1,...,t_s^*$$
$$B_{t_s^*+1} x_{t_s^*} + A_{t_s^*+1} x_{t_s^*+1}^u \le b_{t_s^*+1} \quad (11)$$
$$G_s x_s \le \hat{g}_s$$
$$F_s x_s \le f_s$$

Since the forecast for the future time periods $\mathcal{S}''$ is assumed to remain unchanged, $x_{t_s^*+1}^u$ can be considered as the best available dispatch plan in the time period $t_s^* + 1$ that leads to the optimal schedule in the future $\mathcal{S}''$. The optimal solution $x_s^*$ to $SP_s$ can be thought of as an optimal outcome in the current sub time horizon given the past and future boundary conditions.

The proposed RT pricing problem is defined as follows:

$$PP_s(\pi_{t_s^0-1:t_s^0}^u, \pi_{t_s^*:t_s^*+1}^u) := \min_{x_s} c_s^\mathsf{T} x_s - \pi_{t_s^0-1:t_s^0}^{u,\mathsf{T}} A_{t_s^0} x_{t_s^0} - \pi_{t_s^*:t_s^*+1}^{u,\mathsf{T}} B_{t_s^*+1} x_{t_s^*}$$
$$\text{s.t. } G_s x_s \le \hat{g}_s \quad (\lambda_s^*)$$
$$F_s x_s \le f_s \quad (12)$$
$$B_t x_{t-1} + A_t x_t \le b_t, \quad t = t_s^0 + 1,...,t_s^*$$

where the offer adjustments $-A_{t_s^0}^\mathsf{T} \pi_{t_s^0-1:t_s^0}^u$ and $-B_{t_s^*+1}^\mathsf{T} \pi_{t_s^*:t_s^*+1}^u$ are obtained from the shadow prices of the intertemporal constraints in the forward market. The market clearing price LMP in the current sub time horizon $\mathcal{S}$ can be constructed from the shadow prices of the system level constraint in problem (12) as follows:

$$LMP_s^* = G_s^\mathsf{T} \lambda_s^* \quad (13)$$

*2) Properties of the RTM under the perfect forecast*

In this subsection, we discuss the properties of the RTM clearing results if the forecast is accurate.

**Property 3:** Under the perfect forecast condition, the RT clearing quantities are identical to the forward clearing quantities.

Proof: We prove it by contradiction. Suppose that the forward schedule $x_s^u$ is a feasible but not an optimal solution to the RT scheduling problem $SP_s$. Then, there exists a feasible solution $\tilde{x}_s (\ne x_s^u)$ to $SP_s$ such that $c_s^\mathsf{T} \tilde{x}_s < c_s^\mathsf{T} x_s^u$. Under the perfect forecast assumption, the dispatch solution $(x_{s'}^u, \tilde{x}_s, x_{s''}^u)$ is a feasible solution to the forward problem (1)-(4), and
$$c_{s'}^\mathsf{T} x_{s'}^u + c_s^\mathsf{T} \tilde{x}_s + c_{s''}^\mathsf{T} x_{s''}^u < c_{s'}^\mathsf{T} x_{s'}^u + c_s^\mathsf{T} x_s^u + c_{s''}^\mathsf{T} x_{s''}^u,$$
which contradicts the fact that $(x_{s'}^u, x_s^u, x_{s''}^u)$ is the optimal solution to the forward problem (1)-(4). Therefore, $x_s^u$ is an optimal solution to $SP_s$.

Now we prove the converse. Suppose that the RT schedule $x_s^*$ is a part of a feasible not optimal solution to the forward problem (1)-(4). Then, it implies that the feasible forward solution $(x_{s'}^u, x_s^*, x_{s''}^u)$ satisfies the following inequality:
$$c_{s'}^\mathsf{T} x_{s'}^u + c_s^\mathsf{T} x_s^* + c_{s''}^\mathsf{T} x_{s''}^u > c_{s'}^\mathsf{T} x_{s'}^u + c_s^\mathsf{T} x_s^u + c_{s''}^\mathsf{T} x_{s''}^u,$$
or equivalently, $c_s^\mathsf{T} x_s^* > c_s^\mathsf{T} x_{s'}^u$, which contradicts the fact that $x_s^*$ is the optimal solution to $SP_s$. Q.E.D.

*Remark:* Property 3 indicates that the RT clearing quantities are an optimal solution to the overall multi-period problem even though the RTM optimizes over a sub time horizon rather than the entire time horizon.

**Property 4**: Under the assumption that the forecast for the future time horizon is accurate, the RT clearing quantity of the current sub time horizon ensures the reliability needs for the future time horizon.

Proof: If the forecast for the future time horizon $\mathcal{S}''$ is accurate, then the forward schedule $x_{s''}^u$ must be a feasible dispatch solution in the future time horizon. By setting one of the boundary conditions of the RT scheduling problem at the forward clearing quantity $x_{t_s^*+1}^u$, it can ensure that the schedule in the current sub time horizon, $x_s^*$, is able to transition to the schedule in the following time period, $x_{t_s^*+1}^u$. Hence, the RT schedule $x_s^*$ guarantees a reliable future solution $x_{s''}^u$. Q.E.D.

**Property 5:** Under the perfect forecast condition, the RT clearing prices are identical to the forward clearing prices, and induce dispatch-following incentives.

Proof: We first show the consistency between the forward and RT clearing prices. The forward clearing price $LMP_t^u := G_t^\mathsf{T} \lambda_t^u$, $t = 1,...,T$ is an optimal solution to the following forward dual problem:

$$\max_{(\pi,\lambda,\mu)\leq 0} \sum_{t\in\mathcal{T}} \left( b_t^\mathsf{T} \pi_{t-1:t} + g_t^\mathsf{T} \lambda_t + f_t^\mathsf{T} \mu_t \right) \tag{14}$$

$$\text{s.t.} \quad A_t^\mathsf{T} \pi_{t-1:t} + B_{t+1}^\mathsf{T} \pi_{t:t+1} + G_t^\mathsf{T} \lambda_t + F_t^\mathsf{T} \mu_t = c_t \tag{15}$$
$$t = 1,\ldots,T-1$$

$$A_T^\mathsf{T} \pi_{T-1:T} + G_T^\mathsf{T} \lambda_T + F_T^\mathsf{T} \mu_T = c_T \tag{16}$$

Under the perfect forecast condition, the RT clearing prices $LMP_s^* := G_s^\mathsf{T} \lambda_s^*$ in the sub time horizon $\mathcal{S}$ can be obtained by solving the following dual problem of $PP_s$, which uses the forward dual solution $\pi_{t_s^0-1:t_s^0}^u$ and $\pi_{t_s^*:t_s^*+1}^u$ as the boundary conditions:

$$\max_{(\pi_s,\lambda_s,\mu_s)\leq 0} \sum_{t=t_s^0+1}^{t_s^*} \left( b_t^\mathsf{T} \pi_{t-1:t} + g_t^\mathsf{T} \lambda_t + f_t^\mathsf{T} \mu_t \right) \tag{17}$$
$$+ g_{t_s^0}^\mathsf{T} \lambda_{t_s^0} + f_{t_s^0}^\mathsf{T} \mu_{t_s^0}$$

$$\text{s.t.} \quad A_t^\mathsf{T} \pi_{t-1:t} + B_{t+1}^\mathsf{T} \pi_{t:t+1} + G_t^\mathsf{T} \lambda_t + F_t^\mathsf{T} \mu_t = c_t \tag{18}$$
$$t = t_s^0+1,\ldots,t_s^*-1$$

$$A_{t_s^0}^\mathsf{T} \pi_{t_s^0-1:t_s^0}^u + B_{t_s^0+1}^\mathsf{T} \pi_{t_s^0:t_s^0+1} + G_{t_s^0}^\mathsf{T} \lambda_{t_s^0} + F_{t_s^0}^\mathsf{T} \mu_{t_s^0} = c_{t_s^0} \tag{19}$$

$$A_{t_s^*}^\mathsf{T} \pi_{t_s^*-1:t_s^*} + B_{t_s^*+1}^\mathsf{T} \pi_{t_s^*:t_s^*+1}^u + G_{t_s^*}^\mathsf{T} \lambda_{t_s^*} + F_{t_s^*}^\mathsf{T} \mu_{t_s^*} = c_{t_s^*} \tag{20}$$

Similar to the contradiction proof used in Property 3, the equivalence between $\lambda_s^*$ and $\lambda_s^u$ can be obtained by showing that $\lambda_s^*$ must be the optimal solution to the forward dual problem (14)-(16), and $\lambda_s^u$ is also the optimal solution to the RT pricing dual problem (17)-(20). The details are omitted here to avoid repetitiveness. As a result, it implies that the forward clearing price $LMP_s^u$ is the same as $LMP_s^*$ in the sub time horizon.

Furthermore, because of the equivalency between the forward and RT clearing prices as well as quantities, it is easy to show that the solution to the resource profit maximization problem (5)-(7) matches the solution to the scheduling problem $SP_s$. This implies that each resource in the RTM has incentive to follow the schedule path determined in the forward market. Q.E.D.

*Remark:* The offer adjustments $-A_{t_s^0}^\mathsf{T} \pi_{t_s^0-1:t_s^0}^u$ and $-B_{t_s^*+1}^\mathsf{T} \pi_{t_s^*:t_s^*+1}^u$ can be considered as the boundary conditions for the pricing problem $PP_s$. Under such boundary conditions, the ISO is able to reproduce the full multi-period prices by solving the pricing problem in a shorter time horizon under the perfect forecast condition.

*3) Impact of uncertainty on the opportunity cost estimation and system feasibility*

When the forecast of the uncertain parameter $g$ is imperfect, the properties in the previous subsection do not necessarily hold. To obtain the actual intertemporal opportunity costs, the ISO would have to frequently re-optimize the multi-period problem with the updated parameter $g$ in the RTM, which is impractical and computationally challenging for a large system. In practice, if the realization of the parameter $g$ does not deviate significantly from the forecasted value, the offer adjustments $-A_{t_s^0}^\mathsf{T} \pi_{t_s^0-1:t_s^0}^u$ and $-B_{t_s^*+1}^\mathsf{T} \pi_{t_s^*:t_s^*+1}^u$ can be considered as the estimations of the intertemporal opportunity costs as indicated in Property 6 below. If the realization deviates significantly from the forecast value, then the estimations become inaccurate, re-optimization should be performed.

Before presenting Property 6, we introduce two terms. The *backward-looking intertemporal opportunity cost* is how much the overall profit of all the resources in the past would have changed if every resource is allowed to deviate its production at the start of current look-ahead time horizon by one more unit. The *forward-looking intertemporal opportunity cost* is the overall profit change in the future if every resource is allowed to deviate its production at the end of current look-ahead time horizon by one more unit. Opportunity cost is typically forward looking in that it measures the value that a resource sacrifices at the time the decision is made and beyond. We want to provide justifications for the backward-looking intertemporal opportunity cost here. At any moment in time, a resource can look either forward or backward. One looks backward in time in a perspective of the alternative of leaving the market. One looks forward in time in a perspective of the alternative decisions to be made and the consequences of each. The goal of the market is to incentivize resources to follow the scheduling path for the entire time horizon desired by the ISO. This means that the ISO needs to keep track of the impact of the dispatch made in the past on the present choice. Otherwise, resources would wish they could renege on their previous dispatch actions. By incorporating both the backward-looking and forward-looking intertemporal opportunity costs in the RT pricing model, resources have incentives to follow the entire scheduling path desired by the ISO.

**Property 6:** The offer adjustments $-A_{t_s^0}^\mathsf{T} \pi_{t_s^0-1:t_s^0}^u$ and $-B_{t_s^*+1}^\mathsf{T} \pi_{t_s^*:t_s^*+1}^u$ are approximations of the backward-looking and forward-looking intertemporal opportunity costs, respectively.

Proof: If every resource were able to re-optimize in the past time horizon $\mathcal{S}'$ given that the schedule at the boundary time period $t_s^0$ is set to be $x_{i,t_s^0}$, then it would solve the following profit maximization problem $U_i(x_{i,t_s^0})$, which is parameterized by $x_{i,t_s^0}$:

$$U_i(x_{i,t_s^0}) := \max_{x_{i,s'}} \left( LMP_{s'}^* - c_{i,s'} \right)^\mathsf{T} x_{i,s'} \tag{21}$$

$$\text{s.t.} \quad B_{i,t} x_{i,t-1} + A_{i,t} x_{i,t} \leq b_{i,t},\ t=1,\ldots,t_s^0 \tag{22}$$

$$G_{i,s'} x_{i,s'} \leq g_{i,s'}^*,\ F_{i,s'} x_{i,s'} \leq f_{i,s'} \tag{23}$$

The overall profit of all the resources in the past time horizon is

$$\sum_i U_i(x_{i,t_s^0}) = \max_{x_{s'}} \sum_i -c_{i,s'}^\mathsf{T} x_{i,s'} + LMP_{s'}^* \sum_i x_{i,s'},$$

which is equivalent to the following production cost minimization problem in the past time horizon given the fact that supply and demand are balanced, i.e. $\sum_i x_{i,s'} = 0$:

$$\hat{U}(x_{t_s^0}) := \min_{x_{s'}} c_{s'}^{\mathsf{T}} x_{s'} \tag{24}$$

s.t. $B_t x_{t-1} + A_t x_t \leq b_t,\ t = 1,...,t_s^0$ (25)

$G_{s'} x_{s'} \leq g_{s'}^*,\ F_{s'} x_{s'} \leq f_{s'}$ (26)

It follows by the definition that the backward-looking intertemporal opportunity cost is equal to the subgradient[4] of the function $\sum_i U_i(x_{i,t_s^0})$, which is not necessarily differentiable. Because $\sum_i U_i(x_{i,t_s^0}) = \hat{U}(x_{t_s^0})$, we have the backward-looking intertemporal opportunity cost is equivalent to the subgradient of $\hat{U}(x_{t_s^0})$, i.e. $\partial \hat{U}(x_{t_s^0})$.

It is trivial to show that the forward optimal dual solution ($\pi_{0:1}^u,...,\pi_{t_s^0-1:t_s^0}^u, \lambda_{s'}^u, \mu_{s'}^u$) is dual feasible to $\hat{U}(x_{t_s^0})$. By Proposition 2.1 in [26], the linear function

$$D_{s'}(x_{t_s^0}) := \left(b_{s'}^{\mathsf{T}} \pi_{s'}^u + g_{s'}^{*,\mathsf{T}} \lambda_{s'}^u + f_{s'}^{\mathsf{T}} \mu_{s'}^u\right) - \pi_{t_s^0-1:t_s^0}^{u,\mathsf{T}} A_{t_s^0} x_{t_s^0} \tag{27}$$

is a lower cut for $\hat{U}(x_{t_s^0})$ at $x_{t_s^0}^u$, i.e. for every $x_{t_s^0}$, we have $\hat{U}(x_{t_s^0}) \geq D_{s'}(x_{t_s^0})$ and the distance between the values of $\hat{U}(x_{t_s^0})$ and of the cut $D_{s'}(x_{t_s^0})$ at $x_{t_s^0}^u$ is at most $\epsilon$, where $\epsilon = \hat{U}(x_{t_s^0}^u) - D_{s'}(x_{t_s^0}^u)$. Thus,

$$\hat{U}(x_{t_s^0}) - \hat{U}(x_{t_s^0}^u) \geq -\pi_{t_s^0-1:t_s^0}^{u,\mathsf{T}} A_{t_s^0}(x_{t_s^0} - x_{t_s^0}^u) - \epsilon,$$

which implies that $-\pi_{t_s^0-1:t_s^0}^{u,\mathsf{T}} A_{t_s^0}$ is an $\epsilon$-subgradient[5] of $\hat{U}(x_{t_s^0})$ at $x_{t_s^0}^u$. This means that $-A_{t_s^0}^{\mathsf{T}} \pi_{t_s^0-1:t_s^0}^u$ is an approximation of $\partial \hat{U}(x_{t_s^0})$. In other words, $-A_{t_s^0}^{\mathsf{T}} \pi_{t_s^0-1:t_s^0}^u$ is an approximation of the backward-looking intertemporal opportunity cost.

The optimal future profit for each resource can be defined by $V_i(x_{i,t_s^*})$, which is parameterized by $x_{i,t_s^*}$:

$$V_i(x_{i,t_s^*}) := \max_{x_{i,s''}} \left(LMP_{s''}^* - c_{s''}\right)^{\mathsf{T}} x_{i,s''} \tag{28}$$

s.t. $B_{i,t} x_{i,t-1} + A_{i,t} x_{i,t} \leq b_{i,t},$
$t = t_s^* + 1,...,T$ (29)

$G_{i,s''} x_{i,s''} \leq g_{i,s''}^*,\ F_{i,s''} x_{i,s''} \leq f_{i,s''}$ (30)

Similarly, the overall profit of all the resources in the future time horizon is $\sum_i V_i(x_{i,t_s^*})$, which is equivalent to the following production cost minimization problem in the future time horizon:

$$\hat{V}(x_{t_s^*}) := \min_{x_{s''}} c_{s''}^{\mathsf{T}} x_{s''} \tag{31}$$

s.t. $B_t x_{t-1} + A_t x_t \leq b_t,\ t = t_s^* + 1,...,T$ (32)

$G_{s''} x_{s''} \leq g_{s''}^*,\ F_{s''} x_{s''} \leq f_{s''}$ (33)

It follows by definition that the forward-looking intertemporal opportunity cost is equal to $\partial \hat{V}(x_{t_s^*})$. Similar to the proof of the first half of the property, it can be shown that $-B_{t_s^*+1}^{\mathsf{T}} \pi_{t_s^*:t_s^*+1}^u$ is an approximation of $\partial \hat{V}(x_{t_s^*})$. Q.E.D.

*Remark:* Property 6 suggests that the offer adjustments $-\pi_{t_s^*:t_s^*+1}^{u,\mathsf{T}} B_{t_s^*+1} x_{t_s^*}$ and $-\pi_{t_s^0-1:t_s^0}^{u,\mathsf{T}} A_{t_s^0} x_{t_s^0}$ supplements the RT pricing problem with the information associated with the past and future time periods so that the pricing problem does not need to include the detailed constraints related to the past and future.

Next, let us focus on the RT schedules. If the updated forecast $\hat{g}_s$ in the current sub time horizon is significantly different from the forward forecast, then the RT scheduling problem $SP_s$ may not be feasible under the given boundary conditions. Similarly, if the forecasts for the future time horizon deviate significantly from the ones in the forward market, then the forward schedule for the remaining time horizon may not serve as a reliable guideline for the RT processes and needs to be adjusted. One possible approach is to re-optimize the remaining time horizon $\mathcal{S} \cup \mathcal{S}''$ based on the latest information. The adjusted schedules and intertemporal opportunity costs would be used as the new guideline for the scheduling and pricing processes in the remainder of the time horizon.

### C. The multi-settlement system

In general, a multi-settlement system is the settlement of a sequence of markets for a product that includes at least one forward market. In the forward market, buyers and sellers may conclude financial contracts for later delivery in a RTM, which involves the actual delivery of the product.

Under the proposed framework, the clearing quantity of each time period $t$ at the forward market is settled at the corresponding price: $LMP_t^u \times x_t^u$. In a RTM with a sub time horizon $\mathcal{S}$, the deviation from the previous market clearing quantity, $x_t^{pre}$, is settled at the RT clearing price, i.e. $LMP_{t(s)}^* \times \left(x_{t(s)}^* - x_t^{pre}\right)$, for $t(s) \in \mathcal{S}$. The previous market can be either the forward market or the RTM preceding the current sub time horizon $\mathcal{S}$.

---

[4] If the function $\sum_i U_i(x_{i,t_s^0})$ is differentiable, then the backward-looking intertemporal opportunity cost is equal to the derivative of the function $\sum_i U_i(x_{i,t_s^0})$ by its definition. However, the function $\sum_i U_i(x_{i,t_s^0})$ is not necessarily differentiable, then the subgradient, which is a generalization of the derivative of a function, represents the backward-looking intertemporal opportunity cost. The vector $v$ is a subgradient of a function $f$ at $x$ if $f(x) + v^T(y-x) \leq f(y)$ for all $y$. The set of all subgradients can be written as $\partial f(x)$.

[5] A vector $\omega$ is a $\epsilon$-subgradient of a convex function $f$ at the point $x$ if for every $y$, it satisfies the following inequality: $f(x) + \omega^T(y-x) - \epsilon \leq f(y)$.
7

A benefit of the multi-settlement system is that it reduces risk exposure for market participants. The forward schedule is financially binding. The RT price is used only to settle the deviation from the previous market clearing. By locking in the forward clearing prices, market participants are only exposed to the RT price volatility for the deviation quantity.

## IV. NUMERICAL EXPERIMENTS

In the numerical experiments, we first consider a simple one-node system with an ESR to illustrate the proposed market clearing process and the benefits of coordinating the DAM and RTM. Next, we use the real-world ISO NE system to compare the proposed design with other existing approaches in terms of computational efficiency, economic efficiency, load-following incentives and reliability.

*Case A: One-node system under the perfect as well as imperfect forecast*

In the one-node system, three conventional generators and an ESR are dispatched to serve the load. The characteristics of each generator and the ESR are summarized in Table 1. The charging/discharging efficiency of the ESR is assumed to be 1 for the sake of simplicity. The ESR can seamlessly go from charging to discharging.

|  | Gen1 | Gen2 | Gen3 | ESR |
|---|---|---|---|---|
| Offer Price ($/MWh) | 10 | 63 | 100 | 9 |
| Bid Price ($/MWh) | - | - | - | 5 |
| EcoMax (MW) | 40 | 40 | 30 | 15 |
| EcoMin (MW) | 0 | 0 | 0 | -15 |
| Max SOC (MWh) | - | - | - | 12 |
| Initial SOC (MWh) | - | - | - | 6 |

**Table 1: The characteristics of the resources**

The full time horizon consists of 8 time periods. The proposed market design is applied under two conditions, one assuming the perfect forecast and the other assuming the imperfect forecast. The realized RT load under these two conditions is summarized in Table 2. Under the perfect forecast, the RT realized load is the same as the forecast value used in the forward market. There are two peak loads occurring at t5 and t8.

| Time | RT load (MW) | |
|---|---|---|
|  | Perfect forecast | Imperfect forecast |
| $t1$ | 24 | 25.8 |
| $t2$ | 46 | 42.9 |
| $t3$ | 70 | 64.0 |
| $t4$ | 83 | 83.0 |
| $t5$ | 98 | 91.6 |
| $t6$ | 60 | 56.3 |
| $t7$ | 77 | 72.8 |
| $t8$ | 102 | 111.1 |

**Table 2: The realized load under the perfect and imperfect forecasts**

Figure 2 summarizes the clearing prices as well as ESR clearing quantities in the forward and RT markets under the perfect forecast condition. The coordinated scheduling process is shown at the left side of Figure 2. In the forward market, the ESR is scheduled to charge to its max state-of-charge (SOC) capacity at t1, and remains idle until the first peak load at t5, when it is scheduled to fully discharged. Then, the ESR is charged at t6 and t7, when the load is low. Finally, the ESR is fully discharged at the last time period since any remaining SOC would be worthless afterwards. In the RTM, a sub time horizon contains three time periods, and the actual load is the same as the forecast load under the perfect forecast condition. Each RT scheduling problem utilizes the forward schedules to guide the dispatch. The resulting schedule in each RT scheduling problem matches the corresponding forward schedule, as claimed in Property 3.

The right side of Figure 2. illustrates the coordinated pricing processes. The forward clearing prices and shadow prices of the SOC constraints are obtained by using the dual optimal solutions of the forward scheduling problem. In the RTM, each RT pricing problem takes into account the shadow prices of the SOC constraints or the intertemporal opportunity cost of the ESR. The resulting RT clearing prices are consistent with the forward prices, as proved in Property 5.

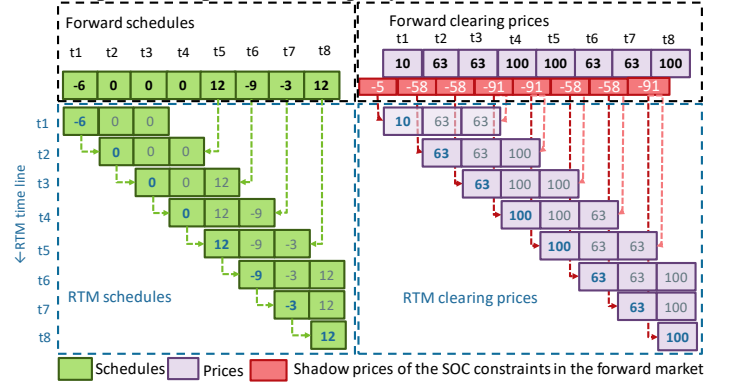

**Figure 2: The ESR schedules and clearing prices under the perfect forecast**

Under the imperfect forecast condition, the RT market clearing results for the ESR are shown in Figure 3. Because the RT load deviates from the forecast, the RT schedule of ESR is different from the forward schedule at t4 and t5 to account for the forecast errors. In other RT time intervals, the schedules of the conventional generators, which is not shown in Figure 3, are adjusted from their forward schedules to meet the actual load. Although the load forecast is not perfect, the resulting RTM clearing prices are identical to the forward clearing prices in this case as shown in Figure 3. In general, the RTM clearing prices are not necessarily equal to the forward prices in the presence of uncertainties.

In order to evaluate the efficiency of the RT schedules and prices under the proposed approach, we compare them with the outcomes under the "perfect information", which assumes that the realized RT load is known to the ISO. The multi-period market clearing model used in the forward market is run with the actual realized RT load to obtain the perfect schedules and prices, which are presented in Figure 3. It can be observed that the RT schedules and RT clearing prices of the ESR are consistent with the outcomes under the "perfect information". The consistency results are also observed for the conventional generators. This implies that the proposed approach results in economic efficient schedules, and all resources have incentives to follow dispatch instructions in this case.

By comparing the forward and perfect schedules in Figure 3, it can be observed that they are very similar. This means that the forward schedules serve as a reliable guideline for the RT dispatches, which contributes to the efficient RT schedules. Similarly, the forward shadow prices of the SOC constraints are identical to the ones under the "perfect information" as shown in Figure 3. This suggests that the forward market accurately captures the intertemporal opportunity cost in this test case, which helps the convergence between the RTM clearing prices and the perfect prices.

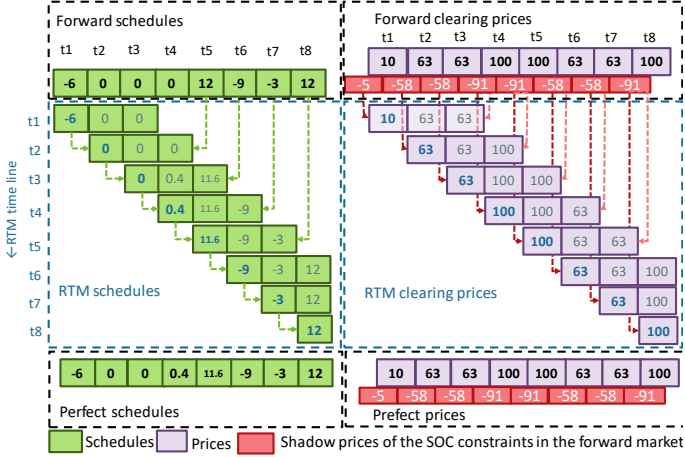

**Figure 3: The ESR schedules and clearing prices under the imperfect forecast**

*Case B: ISO NE system under the imperfect forecast*

In this case, we study the following five approaches: 1) the myopic approach, 2) the first-period only settlement approach [15] with a rolling 2-period RTM, which is denoted by "first only" for short, 3) the proposed approach with rolling 2-period RTM, 4) Hogan's approach [18] and 5) Hua's approach [19] on the real-world ISO NE system. The myopic approach is considered as a baseline for comparing the performance of other alternative approaches. In this study, pumped-storage units with a total of 2000 MW capacity participate in the market. Each pumped-storage unit needs to meet an end-of-day target SOC. Generators are also constrained by their ramping capabilities.

Under each approach, we perform the scheduling and pricing runs to meet 25 randomly realized load scenarios, which uniformly deviate 10% above the forecasted load. We compare the alternative approaches in the following aspects:

(a) Social surplus (SS): SS is a measurement of economic efficiency of a scheduling scheme. The larger SS is, the more economically efficient a dispatch scheme is. SS equals the sum of the conventional generators' producer surplus (PS), consumer surplus (CS) and ESR surplus (ESRS).

(b) Computation time: In practice, the pricing and scheduling models are typically run every 5 minutes in the RTM. Therefore, it is very important that the models must be solved efficiently for a large-scale system in the RTM application.

(c) The number of instances of constraint violation: a measure of reliability of a scheduling scheme. The number of instances of constraint violation are the total number of market intervals when supply and demand is not balanced, or the end-of-day SOC target is not met. The larger the number is, the less reliable the system is under the corresponding dispatch scheme.

(d) Lost opportunity cost (LOC): a measure of dispatch-following incentives. It is the difference between the maximum possible profit a resource can make and the actual profit obtained by following the ISO's dispatch instruction under the given LMPs. In particular, LOC is defined as follows:

$$LOC = \max_{x_i \in \{(6),(7)\}} \sum_{t \in \mathcal{T}} (LMP_t - c_{i,t}) x_{i,t} - \sum_{t \in \mathcal{T}} (LMP_t - c_{i,t}) x_{i,t}^*$$

If the LOC of a resource is equal to zero, it means that the market clearing scheme provides the best dispatch-following incentives. On the other hand, a large LOC implies that the LMP is either too high (or too low), which might lead to the resource generating more (or less) than the dispatch instructions.

Table 3 compares five approaches in terms of the averages of ESRS, CS, PS and SS. Among all the approaches, the myopic approach benefits the conventional producers most at the expenses of consumer surplus and ESR surplus. Since the same scheduling problem is used under Hua's and Hogan's approaches, they result in the same SS. These two approaches outperform all the other approaches in terms of economic efficiency. Under the proposed approach, although the SS is not as high as the Hogan and Hua's approaches, it improves ESRS drastically. Furthermore, SS under the proposed approach is higher than the myopic approach. This implies the dispatch solutions obtained by the proposed approach are more economically efficient than the myopic solutions.

|  | ESRS | CS | PS | SS |
|---|---|---|---|---|
| Myopic (Million $) | 0.031 | 2,220 | 26 | 2,246 |
| First only | +87.5% | +0.7% | -3.6% | +0.7% |
| Proposed | +119.0% | +1.4% | -0.7% | +1.4% |
| Hua | +40.4% | +6.1% | -31.4% | +5.6% |
| Hogan | +40.4% | +6.1% | -31.4% | +5.6% |

**Table 3: Comparison of the average surpluses**

Table 4 reports the average and maximum CPU time for the RTM pricing and scheduling problems. The pricing and scheduling models were modeled in GAMS version 25.02 using CPLEX version 12.8.0.0. All tests were performed on a 2.6-GHz Intel Core i7-6600U CPU with 16 GB RAM, running 64-bit Windows. The wall-clock time limit for CPLEX is set to be 5 minute to reflect RT market clearing on a 5-minute basis. The computation performance of the proposed approach is similar to the first-period only settlement approach. Both are practical for RT implementation. The myopic approach is the most computationally efficient, but it is not as economically efficient compared to the other approaches as shown in Table 3. On the other hand, Hua and Hogan's approaches are significantly slower than the other approaches. Under Hua's approach, the maximum CPU time of the pricing problem is 205.1 seconds, more than 3 minutes. Under Hogan's approach, 49 pricing problems exceed the CPLEX 5-minute time limit among 600 pricing problems in the simulations. This implies that Hua and Hogan's approaches are impractical for the real-time application.

|  | Pricing problem | | Scheduling problem | |
|---|---|---|---|---|
|  | average | max | average | max |
| Myopic | 1.9 | 2.7 | 1.9 | 2.7 |
| First only | 3.8 | 5.3 | 3.8 | 5.3 |
| Proposed | 3.8 | 5.4 | 3.9 | 5.5 |
| Hua | 89.9 | 205.1 | 52.2 | 197.8 |

| Hogan | 154.9 | 300 | 52.2 | 197.8 |

**Table 4: Comparison of the average CPU time (seconds)**

Next, we compare the approaches in terms of reliability. Although supply and demand are balanced in all hours under the myopic and first-period only settlement approaches, pumped-storage resources are unable to meet the end-of-the-day SOC target in all the scenarios under these approaches. This indicates a major reliability disadvantage of these two approaches due to ignoring the reliability needs for future market intervals or short look-ahead horizon. On the other hand, the system is reliable in all the scenarios under all the other approaches.

Table 5 summarizes the mean of LOC under the five approaches. Hua's approach results in the least total LOC whereas the myopic approach leads to the highest total LOC. Although dispatch-following incentive under the proposed approach is not as strong as compared to Hua's approach, it is a significant improvement of the myopic approach and the first-period only settlement approach. The reduction of the total LOC under Hogan's approach is not as promising as the other approaches. This is because several pricing problems under Hogan's approach are not solved to the optimality under the CPLEX time limit as mentioned earlier, the suboptimal prices provide weak dispatch-following incentives.

|  | ESR | Gen | Consumer | Total |
|---|---|---|---|---|
| Myopic | $84,167 | $11,940 | $1,260 | $97,368 |
| First only | -65% | -85% | +3% | -67% |
| Proposed | -98% | -37% | -13% | -90% |
| Hua | -98% | -93% | -31% | -96% |
| Hogan | -99% | +212% | +209% | -57% |

**Table 5: Comparison of mean of LOC**

To understand the impact of the length of the real-time look-ahead horizon under the proposed approach, we vary the look-ahead horizon from 1 to 4 time periods. Table 6 summarizes the averages of SS, total LOC and pricing problem CPU time, relative to the myopic counterparts. It can be observed that the dispatch-following incentives measured by LOC as well as economic efficiency reflected by SS grow as the number of the look-head period increases. The economic efficiency growth slows down after the look-ahead period reaches 2. In addition, although the CPU time of the pricing problem becomes longer with the longer look-ahead period, it is still computationally tractable for RT applications. No constraints are violated under all look-ahead time horizons, so the proposed approach is reliable. Notice that even when the proposed approach solves a single-period problem, it outperforms the myopic counterpart due to the fact that the proposed approach takes into consideration of the future schedules as well as intertemporal opportunity costs.

|  | Number of look-ahead time periods | | | |
|---|---|---|---|---|
|  | 1 | 2 | 3 | 4 |
| SS | +1.3% | +1.4% | +1.4% | +1.4% |
| LOC | -88% | -90% | -93% | -93% |
| CPU time (seconds) | 1.9 | 3.8 | 6.3 | 7.6 |

**Table 6: The averages of SS, LOC and pricing problem CPU time of the proposed approach relative to the myopic approach under different look-ahead horizons.**

## V. CONCLUSIONS

In this paper, a coordinated multi-period market clearing scheme is proposed for both scheduling and pricing of a power system constrained by intertemporal linkages. Compared to the state-of-the-art multi-period approaches, the proposed scheme is innovative for being the first considering the simultaneous pricing and scheduling coordination between a forward market and RTM. The tight coordination allows the proposed scheme to be able to improve computational efficiency, reliability, economic efficiency, and load-following incentives simultaneously whereas the existing methods have to make tradeoffs among these properties. Furthermore, the proposed models consider a general form, rather than being customized for a specific type of intertemporal constraints. These advantages make the proposed scheme practical, generic and flexible to be applied to a wide spectrum of market timeframes as well as to various resource types.

Compared to the existing myopic method, the dispatch instructions generated by the proposed scheme are more economically efficient and reliable. The market clearing prices of the proposed scheme provide better dispatch-following incentives because the intertemporal opportunity cost is taken into consideration in the pricing problem. These advantages of the proposed multi-period scheme are demonstrated in numerical experiments.